\documentclass[tbtags,reqno]{amsart}
\usepackage[latin1]{inputenc}
\usepackage[english]{babel}
\usepackage[T1]{fontenc}
\usepackage{amssymb,amsthm,amsfonts,amscd,latexsym,hyperref}
\providecommand{\msp}{{\kern.03em}}
\newtheorem{prop}{Proposition}
\providecommand{\N}[1]{\lVert u\rVert^{\mathstrut\scriptscriptstyle#1}}
\providecommand{\NU}{\lVert u\rVert^{\mathstrut\scriptscriptstyle\vphantom{1}}}
\newcommand{\ga}[3]{{\Gamma^{#1}{}_{#2#3}}}
\newcommand{\wpp}{{\wp^{(1)}}}
\newcommand{\pp}{{p^{(1)}}}
\newcommand{\vppi}{{\varpi^{(1)}}}
\newcommand{\ppi}{{\pi^{(1)}}}
\newcommand{\TT}{{T^{\displaystyle\ast}(TM)}}
\begin{document}
    \title[Higher order variational origin of the Dixon's system]{Higher order variational origin of the Dixon's system and its relation to the quasi-classical `Zitterbewegung' in General Relativity}
\author[R.~Ya.~Matsyuk]{Roman~Ya.~\sc Matsyuk
}
\address{Institute for Applied Problems in Mechanics and Mathematics\\15~Dudayev~St., L'viv, Ukraine}
\email{matsyuk@lms.lviv.ua}
\thanks{This work was supported by the grant GA\v CR 201/09/0981 of the Czech Science Foundation.}

\subjclass[2000]{Primary 70H50, 83C10; Secondary 53B21, 53B50}


\keywords{Zitterbewegung, relativistic top, spin, covariant Ostrohrads'kyj mechanics, generalised homogeneous Hamiltonian systems}

\begin{abstract}
We show how the Dixon's system of first order equations of motion for the particle with inner dipole structure together with the side Mathisson constraint follows from rather general construction of the `Hamilton system' developed by Weyssenhoff, Rund and Gr\"asser to describe the phase space counterpart of the evolution under the ordinary Euler-Poisson differential equation of the parameter-invariant variational problem with second derivatives. One concrete expression of the `Hamilton function' leads to the General Relativistic form of the fourth order equation of motion known to describe the quasi-classical `quiver' particle in Special Relativity.
The corresponding Lagrange function including velocity and acceleration coincides in the flat space of Special Relativity with the one considered by Bopp in an attempt to give an approximate variational formulation of the motion of self-radiating electron, when expressed in terms of geometric quantities.
\end{abstract}

\maketitle

\section{Introduction}
Consider a quite popular and fairly general Dixon~\cite{matsyuk:Dixon} system of first order ordinary differential equations\footnote{Our definition of the curvature tensor differs in sign from the one adopted in papers~\cite{matsyuk:Dixon,maciuk:B/R}.}
\begin{equation}\label{maciuk:Dixon}
  \left\{\begin{array}{rcl}
    P'{}_\alpha & = & -\displaystyle\dfrac{1}{2}\,R_{\alpha\beta}{}^{\rho\nu}\dot x^\beta S_{\rho\nu} \\[2.5\jot]
  S'{}_{\alpha\beta} & = & \phantom{-}P_\alpha\dot x_\beta-P_\beta\dot x_\alpha
  \end{array}\right.\qquad S_{\alpha\beta}+S_{\beta\alpha}=0\,,
\end{equation}
written in terms of the covariant derivatives, denoted from here on by {\it prime}.

In the theory of General Relativity such equations should hold along the world line of a quasi-classical particle endowed with the inner angular momentum (said `spin') $S^{\alpha\beta}$, responsible for its dipole structure.

Among several additional side conditions needed to make system~(\ref{maciuk:Dixon}) solvable (see~\cite{maciuk:B/R}), we choose to focus on the one preferred by Mathisson~\cite{matsyuk:Mathisson}
\begin{equation}\label{matsyuk:Mathisson}
\dot x^\rho S_{\rho\alpha}=0.
\end{equation}

Now imagine that someone wishes to construct a sort of `Hamilton' picture of the system~(\ref{maciuk:Dixon}) under the imposed constraint~(\ref{matsyuk:Mathisson}). There exists a non-conventional approach to do this along the following  guidelines. First, try to eliminate the variables~$S^{\alpha\beta}$ by means of taking subsequent differential prolongations of~(\ref{maciuk:Dixon},~\ref{matsyuk:Mathisson}). Further, try to find a variational problem with higher derivatives for thus obtained equations, perhaps, under different constraints. Then pass to the corresponding Hamilton--Ostrohrads'kyj counterpart in terms of the generalized momenta. As the last step, compose some geometric quantities~$S^{\alpha\beta}$ from the canonical variables: the momenta and the velocities. If successful, one regains the system~(\ref{maciuk:Dixon}), with the constraint~(\ref{matsyuk:Mathisson}) already satisfied identically.

{\sl We show to the end of this paper that~(\ref{maciuk:Dixon}) follows from a fairly general setting of the second order parameter-invariant variational problem as its `Hamiltonian' counterpart by the appropriate definition of~$S^{\alpha\beta}$.}

In flat space-time of Special Relativity the differential elimination of the variable~$S^{\alpha\beta}$ from~(\ref{maciuk:Dixon},~\ref{matsyuk:Mathisson}) leads to the fourth order equation of motion

\begin{equation}\label{matsyuk:Riewe}
\ddddot x+\left(k^2-\dfrac{m^2}{\sigma^2}\right)\ddot x=0,\qquad (\dot x\cdot\dot x)=1\,,
\end{equation}
where $k^2=(\ddot x\cdot\ddot x)$ is the first integral of~(\ref{matsyuk:Riewe}), and
\begin{equation}\label{sigma}
\sigma_\alpha=\frac{1}{2\NU}\epsilon_{\alpha\beta\rho\nu}u^\beta S^{\rho\nu}\,.
\end{equation}
Equation~(\ref{matsyuk:Riewe}) was shown by Riewe~\cite{matsyuk:Riewe} and Costantelos~\cite{matsyuk:Costantelos} to describe `Zitterbewegung' (quiver) of a quasi-classical particle.

{\sl We show to the end of this paper that~(\ref{matsyuk:Riewe}) occurs as the natural parameterisation of an Euler--Poisson (said \textit{variational}) equation constrained to the manifold~$k=k_0$, of some parameter--invariant variational problem of the second order.}

The above programme for the flat space-time was carried out in two preceding papers~\cite{matsyuk:Symm2001,matsyuk:DGA9}.\footnote{A technical mistake that slipped in the expression for the Hamilton function in paper~\cite{matsyuk:Symm2001} has been corrected in paper~\cite{matsyuk:DGA9}.}

In present paper rather that go all way round the procedure mentioned above, we merely offer a straightforward generalization oh the `Hamiltonian' depiction obtained in~\cite{matsyuk:DGA9} to the case of (pseudo)Riemannian geometry.

\section{The Gr\"asser--Rund--Weyssenhoff canonical equations.}
In the space of the fourth order Ehresmann velocities~$T^4M$ let us stick to the commonly recognized coordinates $x=\{x^\alpha\}\in M$, $u=\dot x=\frac{dx}{d\tau}(0)$, $\dot u=\frac{d^2x}{d\tau^2}(0)$, $\ddot u=\frac{d^3x}{d\tau^3}(0)$, $\dddot u=\frac{d^4x}{d\tau^4}(0)$. A function~$\mathcal L (x, u, \dot u)$
defined on~$T^2M$, constitutes a parameter-invariant variational problem $\delta\int\mathcal L d\tau=0$ if and only if it satisfies the now well known Zermelo conditions:
\begin{subequations}
\renewcommand{\theequation}{\theparentequation.\arabic{equation}}
\begin{gather}
u^\alpha \dfrac{\partial \mathcal L}{\partial \dot u^\alpha} \equiv 0 \label{matsyuk:Z01} \\
u^\alpha \dfrac{\partial \mathcal L}{\partial u^\alpha} +
2\,\dot u^\alpha \dfrac{\partial \mathcal L}{\partial \dot u^\alpha} -
\mathcal L\equiv 0\,.\label{matsyuk:Z02}
\end{gather}
\end{subequations}

We also recall the definition of the Legendre transformation, that is the mapping $Le\!:T^3M\to \TT $ over $TM$ given by
\begin{subequations}\label{matsyuk:3}
\renewcommand{\theequation}{\theparentequation.\arabic{equation}}
\begin{gather}
\wpp=\dfrac{\partial \mathcal L}{\partial \dot u}\,, \label{matsyuk:3.1}
\\
\wp=\dfrac{\partial \mathcal L}{\partial u} - \mathcal D_\tau \wpp\,,\label{matsyuk:3.2}
\end{gather}
\end{subequations}
where
\begin{equation*}
\mathcal D_\tau = u \frac{\partial }{\mathstrut\partial x} +
\dot u \frac{\partial }{\mathstrut\partial u} +
\ddot u \frac{\partial }{\mathstrut\partial \dot u}
\end{equation*}
denotes the operator of total derivative, and the canonical coordinates in $\TT $ are denoted by $x$, $u$, $p$, $\pp$. Applying $\mathcal D_\tau$ to~(\ref{matsyuk:Z01}) immediately gives that in terms of the mixed set of variables $\{\dot u, p, \pp\}$ the Zermelo conditions look like
\begin{subequations}
\renewcommand{\theequation}{\theparentequation.\arabic{equation}}
\begin{gather}
Z_1\overset{\mathrm{def}}=u^\alpha\wpp{}_\alpha = 0 \label{matsyuk:Z11}\\
Z_2\overset{\mathrm{def}}=u^\alpha\wp_\alpha + \dot u^\alpha\wpp{}_\alpha - \mathcal L =0 \,.
\label{matsyuk:Z12}
\end{gather}
\end{subequations}

The standard Liouville form~$\Lambda$ on $\TT$ reads
\begin{equation*}
\Lambda=p.dx+\pp. du\,.
\end{equation*}

The system of the canonical equations developed in the paper of Gr\"asser~\cite{matsyuk:Grasser}, who took as a basis the works of Rund~\cite{matsyuk:Rund} and Weyssenhoff~\cite{matsyuk:Weyssenhoff}, follow from the exterior differential equation
\begin{equation}\label{matsyuk:G-R-W}
Le^{-1}i_Xd\Lambda=-\lambda Le^{-1}d\mathcal H - \mu Le^{-1}dZ_1\,.
\end{equation}
In this equation $Le^{-1}$ denotes the inverse image operation, acting on forms, with respect to the mapping $Le$, and arbitrary functions $\lambda$ and $\mu$ are defined on $T^3M$. If restricted to the first one of the Zermelo conditions~(\ref{matsyuk:Z11}) along the Legendre transformation,
\begin{equation}\label{matsyuk:Z11-Le}
u.\wpp=0\,,
\end{equation}
the exterior differential equation~(\ref{matsyuk:G-R-W}) defines the Legendre transformation itself, along which the function~$\mathcal H$ keeps being constant. It also produces the Euler--Poisson equation of the fourth order that demonstrates the parametric ambivalence to any local transformation of the independent variable~$\tau$.

Let a vector field~$X$ on $\TT$ along some curve $\big(x(\tau),u(\tau),p(\tau),\pp(\tau)\big)$ be its velocity field,
\begin{equation*}
X=\frac{dx^\alpha}{d\tau}\frac{\partial}{\partial x^\alpha}+\frac{du^\alpha}{d\tau}\frac{\partial}{\partial u^\alpha} + \frac{dp_\alpha}{d\tau}\frac{\partial}{\partial p_\alpha} + \frac{d\pp{}_\alpha}{d\tau}\frac{\partial}{\partial\pp{}_\alpha}\,.
\end{equation*}
Then in the coordinate expression the exterior differential equation~(\ref{matsyuk:G-R-W}) amounts to the following system of the first order differential equations~\cite{matsyuk:Grasser}:
\begin{subequations}\label{matsyuk:8}
\renewcommand{\theequation}{\theparentequation.\arabic{equation}}
\begin{align}
\frac{dx}{d\tau}&=\phantom{-}\lambda\,\frac{\partial \mathcal H}{\partial p}\circ Le \label{matsyuk:8.1}\\
\frac{du}{d\tau}&=\phantom{-}\lambda\,\frac{\partial \mathcal H}{\partial \pp}\circ Le + \mu u \label{matsyuk:8.2}\\
\frac{dp}{d\tau}\circ Le &=-\lambda\,\frac{\partial \mathcal H}{\partial x}\circ Le \label{matsyuk:8.3}\\
\frac{d\pp}{d\tau}\circ Le &=-\lambda\,\frac{\partial \mathcal H }{\partial u}\circ Le - \mu \wpp \,. \label{matsyuk:8.4}
\end{align}
\end{subequations}

\section{Implementing the covariant derivation.}
Let us introduce the following change of local coordinates in~$T^2M$:
\begin{equation*}
\{x,u,\dot u\}\xrightarrow{\phi}\{x,u,u'\}\,,
\end{equation*}
where
\begin{equation}\label{phi}
u'\msp^\alpha=\frac{du^\alpha}{d\tau}+\ga\alpha\beta\rho u^\beta u^\rho\,.
\end{equation}
Inspired by~(\ref{matsyuk:3}), we also may consider the `covariant momenta' given by

\begin{subequations}\label{CovP1P}
\renewcommand{\theequation}{\theparentequation.\arabic{equation}}
\begin{gather}
 \vppi =\dfrac{\partial (\mathcal L\circ\phi^{-1})}{\partial u'}\,,\label{CovP1} \\
 \varpi =\dfrac{\partial (\mathcal L\circ\phi^{-1})}{\partial u}- \vppi{}'\,.\label{CovP}
\end{gather}
\end{subequations}
In~(\ref{CovP1}) we have
\begin{equation}\label{CovP1=P1}
\dfrac{\partial(\mathcal L\circ\phi^{-1})}{\partial u'}=\dfrac{\partial \mathcal L}{\partial \dot u}\circ\phi^{-1}
\end{equation}
by virtue of (\ref{phi}). In~(\ref{CovP}) again on the strength of~(\ref{phi}) one computes
\begin{equation}\label{dL_phi_du}
\dfrac{\partial (\mathcal L\circ\phi^{-1})}{\partial u^\alpha}=\dfrac{\partial \mathcal L}{\partial u^\alpha}\circ\phi^{-1}
-2\ga\rho\alpha\beta u^\beta\dfrac{\partial \mathcal L}{\partial \dot u^\rho}\circ\phi^{-1}.
\end{equation}
On the other hand, the rule for the covariant derivative of a covariant vector says:
\begin{equation}\label{CovP1'}
\vppi\msp'{}_\alpha=\frac{d\vppi{}_\alpha}{d\tau} - \ga\rho\alpha\beta \vppi{}_\rho u^\beta.
\end{equation}
Relation~(\ref{CovP1=P1}) should be understood in terms of the notation~(\ref{matsyuk:3.1}) as
\begin{equation}\label{P1_phi}
\vppi=\mathit{Id}\circ\wpp\circ\phi^{-1},
\end{equation}
from where it immediately follows that also
\begin{equation}\label{dotP1_phi}
\dfrac{\partial \vppi}{\partial \tau}=\dfrac{\partial \wpp}{\partial \tau}\circ\phi^{-1}.
\end{equation}
Inserting~(\ref{P1_phi}) and~(\ref{dotP1_phi}) into~(\ref{CovP1'}) and then together with (\ref{dL_phi_du}) and (\ref{matsyuk:3.1}) into~(\ref{CovP}) gives
\begin{align*}
\varpi & = \left(\dfrac{\partial\mathcal L}{\partial u^\alpha}-\frac{d}{d\tau}\wpp{}_\alpha\right)\circ\phi^{-1} - \ga\rho\alpha\beta u^\beta\wpp{}_\rho\circ\phi^{-1} \\
 & = \wp\circ\phi^{-1} - \ga\rho\alpha\beta u^\beta\wpp{}_\rho\circ\phi^{-1}
\end{align*}
by the definition~(\ref{matsyuk:3.2}).
This suggests the corresponding change of coordinates in the manifold~$\TT$ over~$TM$:
\begin{gather}
\{x,u,p,\pp\}\xrightarrow{\Phi}\{x,u,\pi,\ppi\}\,,\notag\\[2\jot]
\label{Phi}
\left\{
\begin{aligned}
 \ppi{}_\alpha& = \pp{}_\alpha
 , \\
 \pi_\alpha & = p_\alpha
 -\ga\rho\alpha\beta u^\beta\pp{}_\rho
 .
\end{aligned}
\right.
\end{gather}
Thus the Legendre transformation~$Le$ is represented in the coordinates $\{x,u,u',u''\}$ and $\{x,u,\pi,\ppi\}$ by~(\ref{CovP1P}), which is the local expression for
\begin{equation}\label{tildaLe}
\widetilde{Le}=\Phi\circ Le\circ\phi^{-1}.
\end{equation}
\begin{prop}
Let $\mathfrak H=\mathcal H\circ\Phi^{-1}$ depend on $x$, $u$, $\pi$, $\ppi$ through the invariants
\begin{equation}\label{dependencies}
\gamma=u\cdot u,\quad \psi=\pi. u,\quad \eta=\ppi\cdot\ppi
\end{equation}
only. Then the `Hamilton equations'~(\ref{matsyuk:8}) take the shape
\begin{subequations}\label{8_H}
\renewcommand{\theequation}{\theparentequation.\arabic{equation}}
\begin{align}
\frac{dx}{d\tau}&=u \label{8.1.H}
\\
u'&=2\vppi\left(\frac{\partial \mathfrak
H}{\partial \psi}\right)^{-1}\frac{\partial \mathfrak
H}{\partial \eta}\circ\widetilde{Le} + \tilde\mu u \label{8.2.H}
\\
\pi'{}_\alpha\circ\widetilde{Le}&=-R_{\alpha\beta\rho}{}^\nu u^\rho u^\beta\vppi{}_\nu
 \label{8.3.H}
\\
\ppi{}\msp'{}\widetilde{Le}&=-2\,u\,\left(\frac{\partial \mathfrak
H}{\partial \psi}\right)^{-1}\frac{\partial \mathfrak
H}{\partial \gamma} \,\widetilde{Le} -\varpi- \tilde\mu \vppi\,,  \label{8.4.H}
\end{align}
\end{subequations}
where $\tilde\mu$ and $\mu$ from~(\ref{matsyuk:8}) are related by $\tilde\mu=\mu\circ\phi^{-1}$.
\end{prop}

\noindent\textsl{Proof.} First we compute the derivatives of $\gamma$, $\psi$, $\eta$. As far as $\psi$ is a mere contraction of the covariant vector $\pi$ with the contravarient vector $u$, it contains no metric tensor $g_{\alpha\beta}$; thus
\begin{equation}\label{psi/x}
\dfrac{\partial\psi}{\partial x}=0.
\end{equation}
From the Riemannian geometry we recall the formul{\ae} for the partial derivatives of the metric tensor
\begin{align*}
\dfrac{\partial g_{\alpha\beta}}{\partial x^\nu} & =\phantom{-} g_{\alpha\rho}\ga\rho\nu\beta+g_{\beta\rho}\ga\rho\nu\alpha\,,\\
 \dfrac{\partial g^{\alpha\beta}}{\partial x^\nu} & =- g^{\rho\beta}\ga\alpha\rho\nu+g^{\alpha\rho}\ga\beta\rho\nu\,,
\end{align*}
so that
\begin{equation}\label{gamma/x}
\begin{aligned}
\dfrac{\partial \gamma}{\partial x^\nu} & = \phantom{-} 2\,\ga\rho\beta\nu u^\beta u_\rho,  \\
\dfrac{\partial \eta}{\partial x^\nu} & = - 2\,\ga\rho\beta\nu\ppi\msp^\beta\ppi{}_\rho\,.
\end{aligned}
\end{equation}

From~(\ref{Phi}) and~(\ref{tildaLe}) we have
\begin{subequations}
\renewcommand{\theequation}{\theparentequation.\arabic{equation}}
\begin{align}
\label{H/p1}
\dfrac{\partial \mathcal H}{\partial\pp{}_\alpha }\circ\Phi^{-1} & = \dfrac{\partial \mathfrak H}{\partial\ppi{}_\alpha}-\dfrac{\partial \mathfrak H}{\partial\pi_\beta}\ga\alpha\beta\rho u^\rho \\
\label{H/P}
\dfrac{\partial \mathcal H}{\partial p_\alpha }\circ\Phi^{-1} & = \dfrac{\partial \mathfrak H}{\partial\pi_\alpha} \\
\label{H/u}
\dfrac{\partial \mathcal H}{\partial u^\alpha }\circ\Phi^{-1} & = \dfrac{\partial \mathfrak H}{\partial u^\alpha}- \dfrac{\partial \mathfrak H}{\partial\pi_\rho}\ga\beta\alpha\rho \ppi{}_\beta \\
\label{H/x}
\dfrac{\partial \mathcal H}{\partial x^\alpha }\circ\Phi^{-1} & = \dfrac{\partial \mathfrak H}{\partial x^\alpha}- \dfrac{\partial \mathfrak H}{\partial\pi_\rho}\dfrac{\partial\ga\beta\nu\rho}{\partial x^\alpha}\, u^\nu\ppi{}_\beta \,.
\end{align}
\end{subequations}
Moreover, recalling our assumption about the dependencies of $\mathfrak H$ upon its arguments, from~(\ref{dependencies}) we get
\begin{subequations}
\renewcommand{\theequation}{\theparentequation.\arabic{equation}}
\begin{align}
\label{H/P1}
\dfrac{\partial \mathfrak H}{\partial\ppi{}_\alpha } & = 2\,\dfrac{\partial \mathfrak H}{\partial\eta }\ppi{}^\alpha\,,\\
\label{H/U}
\dfrac{\partial \mathfrak H}{\partial u^\alpha} & = 2\,\dfrac{\partial \mathfrak H}{\partial\gamma}\,u_\alpha+\dfrac{\partial \mathfrak H}{\partial\psi}\pi_\alpha\,,\\
\intertext{and further on, by the use of~(\ref{psi/x}) and~(\ref{gamma/x}),}
\label{H/X}
\dfrac{\partial \mathfrak H}{\partial x^\alpha} & = 2\,\dfrac{\partial \mathfrak H}{\partial\gamma}\ga\beta\rho\alpha u^\rho u_\beta - 2\, \dfrac{\partial \mathfrak H}{\partial\eta}\ga\beta\rho\alpha\ppi\msp^\rho\ppi{}_\beta\,.
\end{align}
\end{subequations}

Again from~(\ref{Phi}), and recalling the rules for the covariant derivatives of vectors~(\ref{phi}) and covectors,
\begin{equation}\label{pi'}
\pi'{}_\alpha=\frac{d\pi_\alpha}{d\tau}-\ga\rho\beta\alpha\pi_\rho u^\beta\,,
\end{equation}
we obtain
\begin{subequations}
\renewcommand{\theequation}{\theparentequation.\arabic{equation}}
\begin{gather}
\frac{d\pp}{d\tau}\circ\Phi^{-1}=\frac{d\ppi}{d\tau} \label{p1/t} \\
\begin{split}
\frac{dp_\alpha}{d\tau}\circ\Phi^{-1} & = \frac{d\pi_\alpha}{d\tau} +\dfrac{\partial\ga\nu\beta\alpha }{\partial x^\rho}u^\rho u^\beta \ppi{}_\nu + \ga\nu\beta\alpha u'\msp^\beta\ppi{}_\nu + \ga\nu\beta\alpha u^\beta\ppi\msp'{}_\nu\\
 &\phantom{+ \ga\nu\beta\alpha u^\beta\ppi\msp'{}_\nu}
 -\ga\nu\xi\alpha\ga\xi\rho\beta u^\beta u^\rho\ppi{}_\nu + \ga\xi\beta\alpha\ga\nu\rho\xi u^\beta u^\rho\ppi{}_\nu\,.
\end{split}
\label{p/t}
\end{gather}
\end{subequations}

In view of~(\ref{H/P}) equation~(\ref{matsyuk:8.1}) now becomes
\begin{equation}\label{8.1.1}
\frac{dx}{d\tau}=\lambda\frac{\partial\mathfrak H}{\partial\pi}\circ\Phi\circ Le.
\end{equation}

As far as
\begin{equation}\label{H/pi}
\dfrac{\partial \mathfrak H}{\partial \pi} = \dfrac{\partial \mathfrak H}{\partial \psi}\,u\,,
\end{equation}
and in view of~(\ref{8.1.1}), equation~(\ref{matsyuk:8.1}) transforms into
\begin{equation*}
\frac{dx}{d\tau} = \lambda u \dfrac{\partial \mathfrak H}{\partial\psi }\circ\Phi\circ Le,
\end{equation*}
so the choice
\begin{equation}\label{lambda}
\lambda=\left(\dfrac{\partial \mathfrak H}{\partial \psi}\right)^{-1}\circ\Phi\circ Le
\end{equation}
seems legal. Thus in the sequel we implement the definition
\begin{equation}\label{dx=u}
\frac{dx}{d\tau}=u\,.
\end{equation}

Let us consider equation~(\ref{matsyuk:8.2}). In it we substitute (\ref{H/p1}) together with (\ref{H/P1}) for $\dfrac{\partial \mathcal H}{\partial\pp}$ and afterwards we put $\lambda\dfrac{\partial \mathfrak H}{\partial\pi}\circ\Phi\circ Le=u$ on the strength of~(\ref{8.1.1}) and of~(\ref{dx=u}).
But in view of~(\ref{lambda}) this exactly produces~(\ref{8.2.H}).

Now let us turn to equation~(\ref{matsyuk:8.4}). We apply formula like~(\ref{pi'}) into (\ref{p1/t}) in order to use it in~(\ref{matsyuk:8.4}) together with (\ref{H/u}) and (\ref{H/U}). And one more time we use~(\ref{8.1.1}) accompanied by~(\ref{dx=u}) and then apply~(\ref{lambda}). This amounts to~(\ref{8.4.H}).

The equation~(\ref{matsyuk:8.3}) is the most interesting. This is the evolution equation. Recall that the left hand side there is given by~(\ref{p/t}), in where $\frac{d\pi}{d\tau}$ should be substituted by~(\ref{pi'}) as usual.

As the next step for~$\ppi\msp'\circ\widetilde{Le}$ in (\ref{matsyuk:8.3}) we substitute~(\ref{8.4.H}) and for $u'$ in there we substitute~(\ref{8.2.H}).

In the right hand side we implement formula~(\ref{H/x}) with subsequent use of~(\ref{H/X}), and, moreover, inserting there $\dfrac{\partial \mathfrak H}{\partial \pi}$ from~(\ref{H/pi}) together with~(\ref{lambda}).

The final step consists in grouping the remaining terms to fit in the well known definition of the Riemannian tensor, which in contraction with~$u$ and~$\ppi$ reads
\begin{multline}
\ppi{}_\nu\dfrac{\partial \ga\nu\alpha\beta}{\partial x^\rho} u^\beta u^\rho - \ppi{}_\nu\dfrac{\partial \ga\nu\rho\beta}{\partial x^\alpha} u^\beta u^\rho  \\
+\ppi{}_\nu\ga\nu\rho\xi \ga\xi\alpha\beta u^\beta u^\rho - \ppi{}_\nu\ga\nu\alpha\xi \ga\xi\rho\beta u^\beta u^\rho \\
= \ppi{}_\nu R_{\alpha\beta\rho}{}^\nu u^\rho u^\beta.
\end{multline}
This completes the proof.

\begin{prop}
Equations (\ref{maciuk:Dixon}, \ref{matsyuk:Mathisson}) follow from~(\ref{8.3.H}) if one introduces
\begin{equation}\label{S=u_pi1}
S_{\alpha\beta}=u_\alpha\ppi{}_\beta - u_\beta\ppi{}_\alpha\,.
\end{equation}
\end{prop}

\noindent\textsl{Proof.}
If $S=u\wedge\ppi$, then, taking $\ppi\msp'$ from~(\ref{8_H}), one computes
\begin{equation*}
    S'=u'\wedge\ppi+u\wedge\ppi\msp'=\mu\,u\wedge\ppi-u\wedge\pi-\mu\,u\wedge\ppi=\pi\wedge  u\,.
\end{equation*}
Thus we reasonably identify $P$ with $\pi$ and this way regain the second equation of~(\ref{maciuk:Dixon}).

Equation~(\ref{matsyuk:Mathisson}) is satisfied automatically.

On the basis of skew-symmetric properties of the curvature tensor one rewrites the equation~(\ref{8.3.H}) as follows
\begin{multline*}
    \pi'{}_\alpha=P'{}_\alpha=-R_{\alpha\beta}{}^{\rho\nu}u^\beta u_\rho\ppi{}_\nu\\
    =-\frac12R_{\alpha\beta}{}^{\rho\nu}u^\beta u_\rho\ppi{}_\nu+\frac12R_{\alpha\beta}{}^{\rho\nu}u^\beta u_\nu\ppi{}_\rho\\
    =-\frac12R_{\alpha\beta}{}^{\rho\nu}u^\beta S_{\rho\nu}\quad\text{by (\ref{S=u_pi1})}\,,
\end{multline*}
and thus regains the first equation of~(\ref{maciuk:Dixon}).

\textsl{Remark.} Up to this point any concretization of the expression for the `Hamilton function'~$\mathfrak H$ does not matter.

\section{The `Zitterbewegung' and electron self radiation in covariant Riemannian framework}
Let at last the `Hamilton function' take the following concrete shape:
\begin{equation*}
\mathfrak H=\pi. u + \frac{\N3}{4}\ppi\cdot\ppi - A\,\NU\,.
\end{equation*}
Equations (\ref{8.2.H}) and (\ref{8.3.H}) look like
\begin{align}
\label{_1}
 u'& = \phantom{-}\frac{1}{2}\gamma^{3/2}\vppi + \tilde\mu \,u \\
\label{_2}
 \ppi\msp'\circ\widetilde{Le}& = - \frac{3}{4}\gamma^{1/2}\,(\eta\circ\widetilde{Le})\,u - \varpi - \tilde\mu\,\vppi + \frac{A}{\sqrt\gamma}\,u\,.
\end{align}
Now let us apply to~(\ref{_1}), along the Legendre transformation, the Zermelo condition~(\ref{matsyuk:Z11-Le}) in terms of the covariant variable~$\vppi$, given by~(\ref{P1_phi}). This defines~$\tilde\mu$:
\begin{equation}\label{mu_ultimate}
\tilde\mu=\frac{u\cdot u'}{\N2}\,.
\end{equation}
To solve the system of equations (\ref{_1}, \ref{_2}) for~$\varpi$, it suffice to differentiate~(\ref{_1}) and to substitute there~(\ref{_2}) for~$\vppi\msp'$  and~(\ref{_1}) for~$\vppi$ again.  This resulting in
\begin{equation}\label{ultim_pi}
\boxed{\varpi=6\frac{u'\cdot u}{\N5}\,u'+ \left(2\frac{u''\cdot u}{\N5} - 5\frac{(u'\cdot u)^2}{\N7} - \frac{u'\cdot u'}{\N5}\right)\,u - 2\frac{u''}{\N3} + \frac{A}{\NU}\,u\,.}
\end{equation}

As far as we remember that the equation~(\ref{8.3.H}) keeps the property of the ambivalence to arbitrary transformations of the evolution parameter~$\tau$, we  absolutely may pass to the natural parameter~$s$ for which $u\cdot u=1$ together with $u'\cdot u=0$ and $u''\cdot u=-\,u'\cdot u'$. Then~(\ref{8.3.H}) with~$\varpi$ (\textit{viz.}~~$\pi$) given by~(\ref{ultim_pi}) finally amounts to the following dynamical equation of motion:
\begin{equation}\label{natural_s}
\dfrac{D}{ds}\left[\Big(-3\,\frac{D^2x}{ds^2}\cdot\frac{D^2x}{ds^2} +A\Big)\, \frac{Dx_\alpha}{ds} -2\,\frac{D^3x_\alpha}{ds^3}\right]= -\ppi{}_\nu\,R_{\alpha\beta\rho}{}^\nu \frac{Dx^\rho}{ds}\frac{Dx^\beta}{ds}\,.
\end{equation}

\begin{prop}
In flat space-time equation~(\ref{natural_s}) on the constraint manifold $k=k_0$ reduces  to the Riewe--Costantelos equation~(\ref{matsyuk:Riewe}) of the quasi-classical `Zitterbewegung' with the frequency $\sqrt{\frac{3}{2}k_0^2 - \frac{A}{2}}$ by putting $A=k_0^2 + 2\frac{m^2}{\sigma^2}$ with $\sigma$ from~(\ref{sigma}).
\end{prop}

\noindent\textsl{Proof.} Take respect for the definition $k=\frac{D^2x}{ds^2}\cdot\frac{D^2x}{ds^2}=k_0{}^2$ in~(\ref{natural_s}).

The Lagrange function for the forth order equation~(\ref{8.3.H}) with~$\varpi$ and~$\vppi$ given by (\ref{ultim_pi}), (\ref{_1}), and~(\ref{mu_ultimate}), is:
\begin{equation}\label{L_ultimate}
\mathcal L_k=(k^2+A)\NU\,.
\end{equation}

In 1946 Bopp~\cite{matsyuk:Bopp1946} in the attempt to give an approximate  variational formulation of the self radiating electron in flat space-time of Special Relativity, introduced a Lagrange function that might have been even at that time expressed in terms of the geometric quantity~$k$ --- the first Frenet curvature of the radiating particle's world line. The Bopp's Lagrangian is nothing but~(\ref{L_ultimate}), seen in pseudo-Euclidian coordinates.

\end{document}